\documentclass{easychair}

 \pdfoutput=1
\usepackage{doc}
\usepackage{makeidx}
\usepackage{subfigure}

\newtheorem{thm}{Theorem}
\newtheorem{lem}[thm]{Lemma}
\newtheorem{cor}{Corollary}
\newtheorem{defn}{Definition}

\begin{document}

%
\title{A Cylindrical Basis Function\\
for Solving Partial Differential Equations on Manifolds}%

%
\author{
E.O Asante-Asamani\inst{1}
\and
    Lei Wang \inst{1}
    \and
   Zeyun Yu \inst{2}
}

\institute{
  University of Wisconsin-Milwaukee,
  Department of Mathematical Sciences\\
     \email{eoa@uwm.edu,wang256@uwm.edu}
\and
     University of Wisconsin-Milwaukee,
     Department of Computer Science\\
   \email{yuz@uwm.edu}
 }


\authorrunning{Asante-Asamani, Wang and Yu}

\titlerunning{Cylindrical Basis Function}

\maketitle

\begin{abstract}
 Numerical solutions of partial differential equations (PDEs) on manifolds continues to generate a lot of interest among scientists in the natural and applied sciences. On the other hand, recent developments of 3D scanning and computer vision technologies have produced a large number of 3D surface models represented as point clouds. Herein, we develop a simple and efficient method for solving PDEs on closed surfaces represented as point clouds. By projecting the radial vector of standard radial basis function(RBF) kernels onto the local tangent plane, we are able to produce a representation of functions that permits the replacement of surface differential operators with their Cartesian equivalent. We demonstrate, numerically, the efficiency of the method in discretizing the Laplace Beltrami operator.
\end{abstract}


\setcounter{tocdepth}{2}
{\small

%
%

\pagestyle{empty}

\section{Introduction}
\label{sect:introduction}
Many applications in the natural and applied sciences require the solution of partial differential equations on manifolds. Such applications arise in areas such as computer graphics\cite{Sun2009,Lai2010,Levy2006}, image processing\cite{Turk1991,Diewald2000,Dorsey2000,Witkin1991}, mathematical physics\cite{Fu2011}, biological systems\cite{Schittkowski2008,Hofer2004}, and fluid dynamics\cite{Misra2011,Myers2004,Myers2002}.A lot more interest, especially in computer graphics, has been generated around solution of PDEs on closed 2D manifolds, as these arise as boundaries of 3D objects. The development of high resolution 3D scanning devices, which capture these surfaces as point clouds, have made numerical methods that can be applied directly to point clouds very attractive. 

A class of numerical methods that have been developed to solve PDEs on closed surfaces, involve the expression of differential operators as projections of their Cartesian equivalents onto local tangent planes via a projection operator ($I-\vec{n}\vec{n}^T$). The resulting operators are then discretized using, for example, finite element methods\cite{Dziuk2007}. Recently, the projection method has been extended to PDE's defined on manifolds represented as point clouds \cite{Flyer2009}. For such methods, functions defined on the manifold are represented using radial basis functions (RBF). The surface differential operators are obtained by applying a projection operator to the RBF discretization of the Cartesian equivalent.

 While the projection methods modify Cartesian differential operators, another class of methods embed the surface PDE into $\Re^3$ so that solutions to the embedded problem when restricted to the surface provide the solution on the surface\cite{Greer2006}. Because these methods result in embeded PDEs posed in $\Re^3$, spatial complications arise when the PDEs have to be solved on a restricted surface domain. The closest point method, developed in \cite{Ruuth2008}, attempts to resolve this problem by extending the functions into $\Re^3$ in a way that makes them constant in the normal direction. This allows the simple replacement of surface differential operator by their Cartesian equivalent. The method however requires a high order interpolation at each time step in order to obtain solutions on the surface.
 
  Recently, the orthogonal gradient method presented in\cite{Piret2012}, was introduced to extend this idea to point clouds. Here,$2N$ additional nodes are introduced during the construction of a distance function, to force the function defined on the surface to be constant in the normal direction. The $2N$ nodes are chosen according to an offset parameter $\delta$ which controls their distance from the surface. The accuracy of the method and condition number of the resulting differential matrix is however sensitive to the choice of $\delta$. The use of $2N$ additional nodes, to enforce derivative constraints, increases the computational complexity of forming the interpolation and differential matrices. 

  In this work we propose a modified RBF kernel, the Cylindrical Basis Function (CBF), which is intrinsically constant in the normal direction at each point of a surface. The modified kernel, when used in the representation of smooth functions defined on a closed surface, allows surface differential operators to be replaced by their Cartesian equivalent without the need to impose additional constraints on the function or have an implicit representation of the surface. We also avoid inherent challenges in performing higher order interpolation, typical of embedding techniques, by discretizing operators directly on the manifold. The proposed method is simple and requires the solution of a much smaller linear system, compared to the orthogonal gradient method.
  
%

\section{The Radial Basis Function}
\label{sect:radial basis functions}
Given function data $\{f_k\}_{k=1}^N$ at the node locations $\{x_k\}_{k=1}^N \in \Re^d$, the RBF interpolant $s(x)$ to the data is given as 
\begin{equation}\label{rbfeqn}
s(x) = \sum^N_{i=1} \lambda_i \phi(\parallel x-x_i\parallel),
\end{equation}

where $\phi(\parallel x-x_i\parallel)$ is the RBF kernel centered at the node $x_i$, and $\lambda_i$ are coefficients chosen to satisfy the interpolation conditions
\begin{equation}\label{intcond}
s(x_i) = f(x_i) \hspace{5mm} i=1 \cdots N,
\end{equation}
which is equivalent to solving the linear system 
\begin{equation}\label{rbfsys}
A\vec{\lambda} = \vec{f}.
\end{equation}

$A$ is the matrix with entries $a_{ij} = \phi(\parallel x_i-x_j\parallel)\hspace{2mm}i,j=1\cdots N$, usually referred to as the interpolation matrix. The norm $\parallel.\parallel$ is taken to be the Euclidean norm. Some common kernels are presented in Table \ref{rbftab}.

\begin{table}[h]
\begin{center}
\begin{tabular}{|l|c|c|}
\hline
Name of RBF & Abbreviation & Definition\\
\hline
Multiquadric & MQ & $\sqrt{1+cr^2}$\\
\hline
Inverse Multiqudric & IMQ & $\frac{1}{\sqrt{1+cr^2}}$\\
\hline
Inverse Quadric & IQ & $\frac{1}{1+cr^2}$ \\
\hline
Gaussian & GA & $e^{-cr^2}$\\
\hline
\end{tabular}
\caption{Common radial basis functions with shape parameter c}\label{rbftab}
\end{center}
\end{table}
\section{RBF Discretization of Cartesian Differential Operators}
\label{sect:cartesian}
Given a point cloud $P=(x_1,x_2,\cdots,x_N)$ and data from a smooth function $\{f_k\}_{k=1}^N$ defined on these points, an RBF interpolation of $f$ satisfies,
\begin{equation}
f(x_i) = \sum_{j=1}^N \lambda_j\phi(||x_i-x_j||) \hspace{5mm} i = 1,2,\cdots,N.
\end{equation}
Let $L$ be a differential operator acting on $f$ and $g(x_i)$ the value of $Lf$ at the point $x_i$ then,
\begin{equation}
g(x_i) = \sum_{j=1}^N \lambda_jL\phi(||x-x_j||)|_{x=x_i} \hspace{5mm} i = 1,2,\cdots,N,
\end{equation}
which defines a linear system and can be represented in matrix form as,
\begin{equation}\label{rbfsys1}
B\vec{\lambda} = \vec{g},
\end{equation}
where the differential matrix $B$ has entries $b_{ij} = L\phi(||x-x_j||)|_{x=x_i} i,j = 1,\cdots N$.

The interpolation matrix A in Eq.\eqref{rbfsys} is non-singular (see \cite{Flyer2009}) and permits the substitution of $\lambda = A^{-1}f$ into Eq.\eqref{rbfsys1} leading to $\vec{g} = BA^{-1}\vec{f}$. The differentiation matrix $BA^{-1}$ gives the RBF discretization of L with respect to the point cloud P.

\section{Construction of Surface Differential Operators}
\label{sect:surface operators}
Let $f$ be a smooth function defined on an arbitrary surface $\Gamma$. The gradient of $f$ expanded in the normal, first and second tangent orthogonal coordinates $\{\vec{n},\vec{t_1},\vec{t_2}\}$ at some point $\vec{x} \in \Gamma$ can be expressed as 
\begin{equation}
\nabla f = \partial_n f \vec{n} + \partial_{t_1} f \vec{t_1} + \partial_{t_2} f \vec{t_2}.
\end{equation}
The surface gradient operator $ \nabla_{\Gamma}$ is the projection of the regular gradient onto the local tangent plane at $\vec{x}$. Thus,
\begin{equation}
\nabla_{\Gamma} f = \partial_{t_1} f \vec{t_1} + \partial_{t_2} f \vec{t_2}.
\end{equation}
Two approaches are typically used in the literature \cite{Flyer2009,Piret2012} to obtain surface operators from regular operators. The first involves projecting the Cartesian gradient operator onto the local tangent plane of the surface at some surface point. The surface gradient becomes,
\begin{align*}
\nabla_{\Gamma} f &=\nabla f - \partial_n f\vec{n}\\
                          &= (I-\vec{n}\vec{n}^T)\nabla f.
\end{align*}
Such methods are classified as projection methods. The second approach involves extending the function f, into $\Re^d$ such that it is constant in the normal direction at each point on the surface. As pointed out by \cite{Ruuth2008} under such conditions the surface gradient and Cartesian gradient agree on the surface. Surface differential operators can then be replaced with the simpler Cartesian operators. The Orthogonal Gradient Method enforces this requirement by extending an RBF approximation of the function outside of the surface $\Gamma$, originally having N points, using $2N$ additional points. This increases the complexity of the resulting linear system from $N$ to $3N$. Also the accuracy of the method is influenced by the choice of an offset parameter $\delta$ which controls the proximity of the $2N$ points to the surface.

In order to avoid the increased complexity of introducing 2N additional points to enforce the null gradient condition we propose a modified kernel which projects all radial vectors in the traditional RBF kernel onto the local tangent plane. We refer to this as the Cylindrical Kernel.  Our modified kernel is intrinsically constant in the normal direction and greatly simplifies the construction of surface differential operators. 

%
\section{Cylindrical Basis Function (CBF)}
\label{sect:cylindrical}
\begin{defn} Given a point cloud $P=(x_1,x_2,\cdots,x_N)$ sampled from a smooth manifold, $\Gamma$ and data from a smooth function $\{f_k\}_{k=1}^N$ defined at these points, a CBF interpolant of f is defined as 
\begin{equation}
c(x) = \sum^N_{j=1} \lambda_j \phi(r_j(x))
\end{equation}
where, 
\begin{equation} \label{CBF}
r_j(x) = \parallel (x-x_j)-\vec{n_j}[(x-x_j)\cdot \vec{n_{j}}]\parallel. 
\end{equation}
$\lambda_j$ and $\phi$ are the interpolation coefficients and kernel respectively, $\vec{n_j}$ is the unit normal vector at the point $\vec{x_j}$ and Eq.\eqref{CBF} is the $\lq\lq$cylindrical distance$\lq\lq$ with respect to the Euclidean norm.
\end{defn}

\begin{lem}\label{lemCBF}
Let $P=(x_1,x_2,\cdots,x_n)$ be a point cloud on a 2D-manifold $\Gamma$ embedded in $\Re^3$ then the CBF kernel satisfies for all $x \in \Gamma$
 \begin{equation}
 \nabla \phi(r_i(x))\cdot \vec{n_{y^i}} = 0 \hspace{3mm} \forall i=1,2,\cdots ,N \text{ and } x \in \Re^3.
 \end{equation}
\end{lem} 

Consider the Laplacian of f on $\Gamma$,
\begin{equation}
\bigtriangleup f = (\partial_nf \vec{n}+\partial_{t_1}f \vec{t_1} + \partial_{t_2}f \vec{t_2})\cdot(\partial_nf \vec{n}+\partial_{t_1}f \vec{t_1} + \partial_{t_2}f \vec{t_2})
\end{equation}
expanding out the operator it is clear that the surface Laplacian, $\bigtriangleup_Sf$ is equivalent to the regular Laplacian $\bigtriangleup f$ if $\partial_nf=\partial_n^2f = 0$. Therefore the following corollary is a consequence of Lemma \ref{lemCBF}

\begin{cor}
On a smooth 2D manifold $\Gamma$ the CBF interpolation, $c(x)$ of a smooth function, $f:\Gamma\rightarrow \Re$ satisfies that,
\begin{align}
\nabla c(x) &= \nabla_{\Gamma} c(x)\\
\bigtriangleup c(x) &= \bigtriangleup_{\Gamma} c(x)
\end{align}
\end{cor}

This implies that surface differential operators can now be discretized by simply discretizing their corresponding Cartesian operators. We illustrate the great simplicity of this approach by computing the Laplacian of the cylindrical kernel. We present the main result as follows;
Given $\phi(r)$, the cylindrical kernel, the surface Laplacian is given as,

\begin{equation}\label{laplacian}
\bigtriangleup \phi(r_{ij})= \phi\prime\prime(||r_{ij}||)(1-(\hat{r_{ij}}\cdot\vec{n_j})^2) + \phi'(||r_{ij}||)\frac{1+(\hat{r_{ij}}\cdot\vec{n_j})^2}{||r_{ij}||}
\end{equation}

where $r_{ij} = (x_i-x_j)-\vec{n_j}[(x_i-x_j)\cdot \vec{n_{j}}$ is the projected radial vector computed from the point $x_i$ to the center of the basis function $x_j$. $\hat{r_{ij}}=\frac{r_{ij}}{||r_{ij}||}$  is the normalized projected vector.
\section{Implementation}
\label{sect:Implementation}

We outline the algorithm for discretizing the Laplace-Beltrami Operator on a 2D-manifold.
\begin{enumerate}
\item \textbf{Obtain the surface normals}. Surface normals can be computed for every point in the point cloud $P$ using Principal Component Analysis(PCA). At any point $p_i \in P$ the goal is to find the coordinate system consisting of the normal and two tangential directions to the local coordinate plane  within a neighborhood of $p_i$. Let $X$ be a matrix whose rows $x_j (j=1 \cdots k)$ represent the k nearest neighbors of $p_i$. Denote by $c$ the mean of all the neighbors i.e 
\[ c = \frac{1}{k} \sum_{j=1}^k x_j\] and let $y_j = x_j-c$ be $j^{th}$ row of the matrix $Y$. It is well know that the eigen vectors of the covariance matrix $YY^T$ form an orthogonal basis which are the principal components \cite{Zhao2012}. The eigen vector corresponding to the smallest eigen value is the desired normal direction to the local plane at $p_i$.
\item Assemble the collocation matrix, \[A=\left[
\begin{matrix}
 \phi(r_{11}) & \phi(r_{12})&\dots & \phi(r_{1n})\\
 \phi(r_{21}) & \phi(r_{22})&\dots & \phi(r_{2n})\\
 \hdotsfor{4}\\
 \phi(r_{n1}) & \phi(r_{n2})&\dots & \phi(r_{nn})
 \end{matrix} \right] \] 
 where \[\phi(r_{ij}):= \phi(||(x_i-x_j)-\vec{n_j}[(x_i-x_j)\cdot \vec{n_{j}}]||) \] corresponds to the basis function centered at the $j^{th}$ node and evaluated at node $i$. \\
\item Assemble the differential matrix,\[ B =\begin{pmatrix} 
\bigtriangleup \phi(r_{11}) & \bigtriangleup \phi(r_{12}) & \dots & \bigtriangleup \phi(r_{1n})\\
\bigtriangleup \phi(r_{21}) & \bigtriangleup \phi(r_{22}) & \dots & \bigtriangleup \phi(r_{2n})\\
\hdotsfor{4}\\
\bigtriangleup \phi(r_{n1}) & \bigtriangleup \phi(r_{n2}) & \dots & \bigtriangleup \phi(r_{nn})
\end{pmatrix}\] 
where $\Delta \phi(r_{ij})$ is as defined in Eq.\eqref{laplacian} \\
\item The discrete LB operator is then given as $LB = BA^{-1}.$
\end{enumerate}
%
%
%
%
\section{Numerical Experiments}
\label{sect:experiments}
\subsection{ Eigen Values of LB Operator}
The eigen values and eigen functions of the LB operator provide intrinsic global information that can be used in characterizing the structure of surfaces\cite{Reuter2005,Reuter2009}. We therefore tested the performance of our method in obtaining the spectra of the LB operator on the unit sphere by solving the eigen value problem
\begin{equation}
\Delta_Su = -\lambda u,
\end{equation}
and compared our first 100 eigen values to the exact values. Our discrete operator was computed using a uniform sampling of the unit sphere with nodes ranging from 258-16386. We used a local implementation of the Gaussian radial basis function by choosing the shape parameter c, such that the function vanished outside the support radius. We performed the same experiment using Finite Element method and compared the convergence of both methods using the Euclidean distance as a measure of errors (Figure \ref{eigen1}). As we refine the nodes we see from Figure \ref{eigen} that our computed eigen values line almost perfectly with the true values. With 4098 nodes our computed eigen values are at a distance of 1.4 from the true values while the FEM values are 11.8 units away. Our results indicate that 4 times the number of nodes used for CBF would be required for FEM in order to achieve similar accuracy. We also observed a decline in the accuracy of our solution at 16386 nodes from 4098 nodes. This we believe is as a result of the significantly higher condition number of the collocation matrix. We had to consistently increase the size of the support domain as the number of nodes increased to achieve optimum results. We show the computational time with increasing numbers of nodes. Computations were performed using Matlab R2013a with an intel core i3-4130 8.40 GHZ processor with 8GB of ram.
\begin{figure}
\centering
  \begin{tabular}{cc}
       {\includegraphics[width=2.5 in]{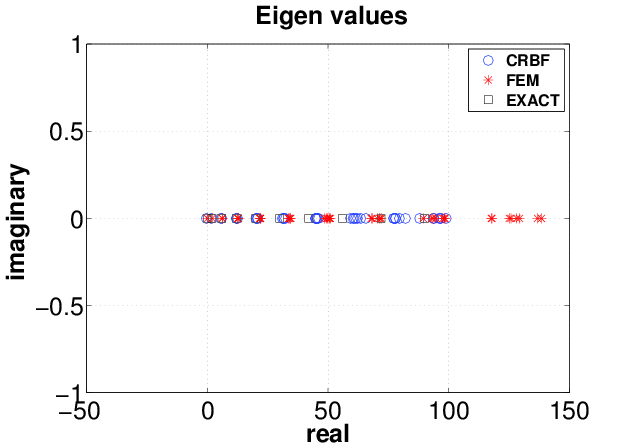}} &  {\includegraphics[width=2.5 in]{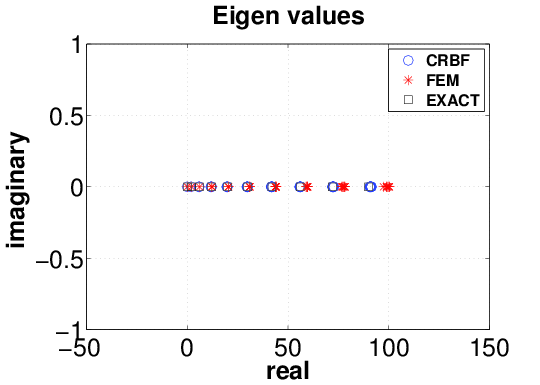}} \\
       {\includegraphics[width=2.5 in]{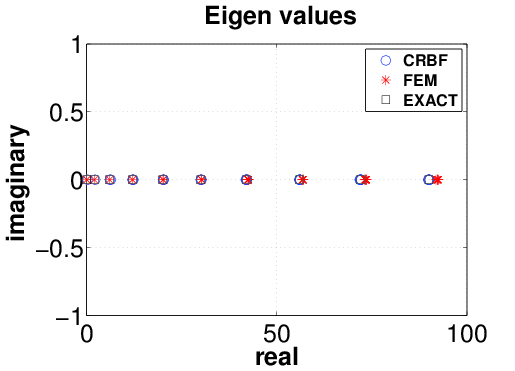}} &  {\includegraphics[width=2.5 in]{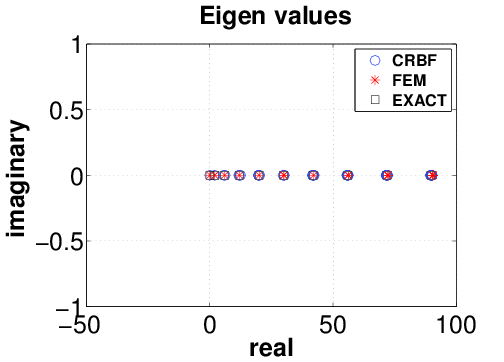}}     
  \end{tabular}
 \caption{first 100 Eigen values of the Laplace Operator on a unit sphere. Number of nodes from left to right 258, 1026, 4098, 16386}\label{eigen}
 \end{figure}
 
 \begin{figure}
\centering
  \begin{tabular}{rcl}
       {\includegraphics[width=2.5 in]{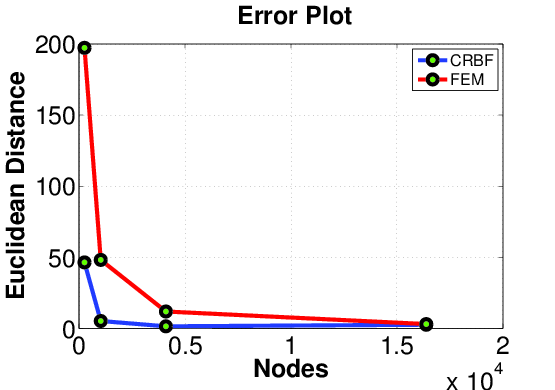}} &  {\includegraphics[width=2.5 in]{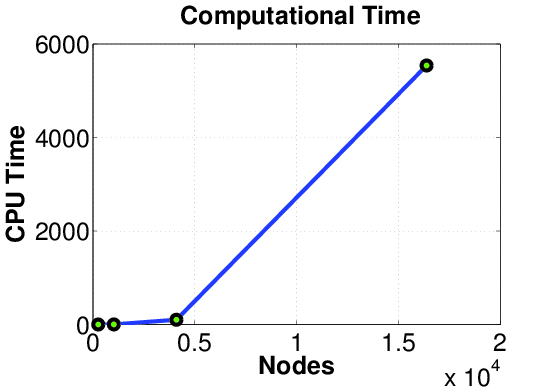}} \\    
  \end{tabular}
 \caption{Superior convergence of CBF and Computational Efficiency. CPU time was measured in secs.}\label{eigen1}
 \end{figure}
\subsection{ Heat Diffusion on Sphere and Molecule}
To verify that our discrete operator accurately captures the heat diffusion of the LB operator we solved the heat equation 
\begin{equation}
u_t = \epsilon \Delta_{\Gamma}u
\end{equation}
on a unit sphere as well as a more complicated molecular surface. We used a Forward Euler time discritization scheme with $\Delta t=0.1h^2$. In Figure \ref{heat1} we show that the diffusion, using our CBF discrete operator, matches well with the exact solution obtained using the spherical harmonics $Y^0_1=\frac{1}{2}\sqrt{\frac{3}{\pi}}x$ across a time duration of 1 sec.


 \begin{figure}
\centering
  \begin{tabular}{c}
       {\includegraphics[width=3.5 in]{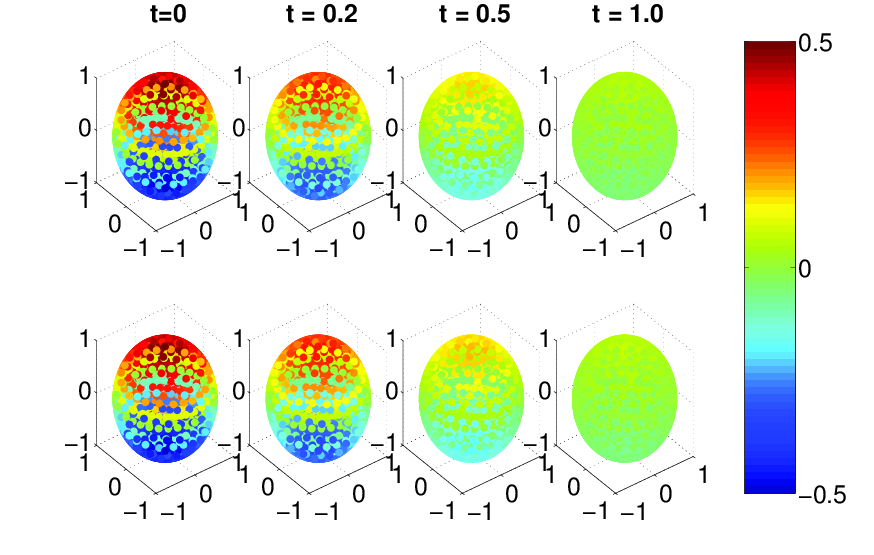}}\\     
  \end{tabular}
 \caption{Heat Evolution on unit sphere with initial condition, the spherical harmonic  $Y_1^0$ and diffusion coefficient $\epsilon=1$ over 1sec duration. (Above) Numerical solution using 4098 nodes. (Below) Exact evolution of sphereical harmonic $Y_1^0$.}\label{heat1}
 \end{figure}
We used total of 7718 point on the molecule surface and allowed the heat to diffuse for 5 secs.  As shown in Figure \ref{heat3}, The heat diffuses as expected from the heat source to neighboring regions. The initial heat distribution was chosen as  a Gaussian bell, $f(x,y,z)=10e^{-4(x-x_1)^2+(y-y_1)^2+(z-z_1)^2}$ centered at a point $(x_1,x_2,x_3)$.

 \begin{figure}
\centering
  \begin{tabular}{cc}
       {\includegraphics[width=2.0 in]{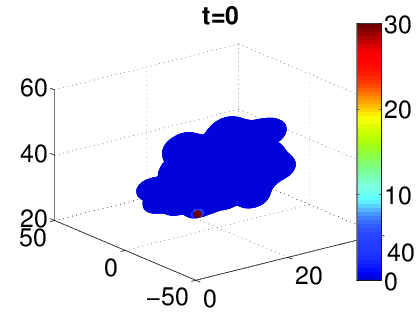}}&  {\includegraphics[width=2.0 in]{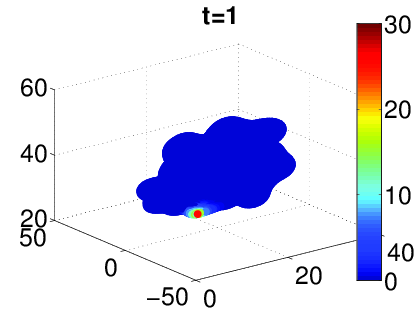}}\\  
       {\includegraphics[width=2.0 in]{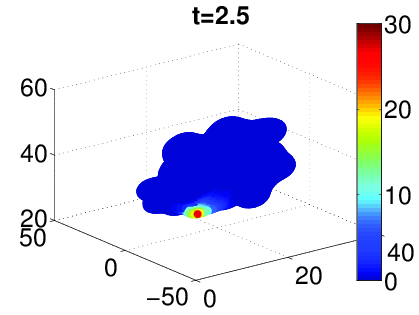}}&  {\includegraphics[width=2.0 in]{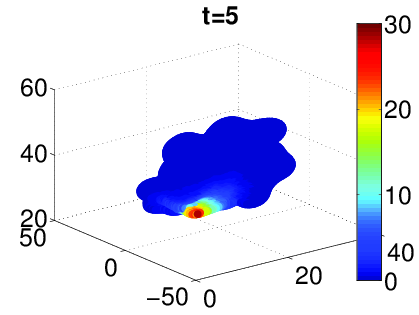}}\\    
  \end{tabular}
 \caption{Heat equation on a molecule with the gaussian bell as initial distribution centered at one of its nodes}\label{heat3}
 \end{figure}


\subsection{ Laplace Smoothing}
Laplace smoothing is a common technique used to smooth surface meshes \cite{Field1988}. Here the Laplace operator acts as a weighted average of the coordinates of each point and its one ring neighbors to draw them towards the barycenter of their common region. To test the applicability of our CBF Laplacian, we used it in smoothing a noisy meshed sphere having 258 nodes. The nodes of the mesh were used as initial condition for the heat equation and updated for 60 iterations. The smoothing equation used was in the form $ u_{i+1} = u_i + \mu \delta LB u_i$ where $ \mu \delta$ acts as the smoothing scale factor with $\mu$ being the diffusion coefficient and $ k$ the time step size. $LB$ is our CBF approximation to the Laplace Beltrami Operator. In Figure \ref{laplacesmooth} we show the smoothing of the triangular mesh on sphere over 60 steps. To measure the effectiveness of our smoothing technique we measured the angle between the point normals at each stage of smoothing to the true normals of the perfectly smooth sphere. Both normals were computed using PCA with 6 neighbours. We choose not to use the position vectors of the point clouds on the perfect sphere as true normals in order to permit a more uniform comparison with the smoothed sphere. We see in Figure \ref{histangles} and Figure \ref{meanangles} that the normals converge to the perfect sphere normals as we iterate to 60 steps. We didn't notice any remarkable improvement in the smoothing after 60 steps. We find the method scales very well with size of mesh nodes as demonstrated in Figure \ref{smoothtime}.
\begin{figure}[ht]
\centering
\subfigure[Rough Mesh]{
    \includegraphics[width=0.35\textwidth]{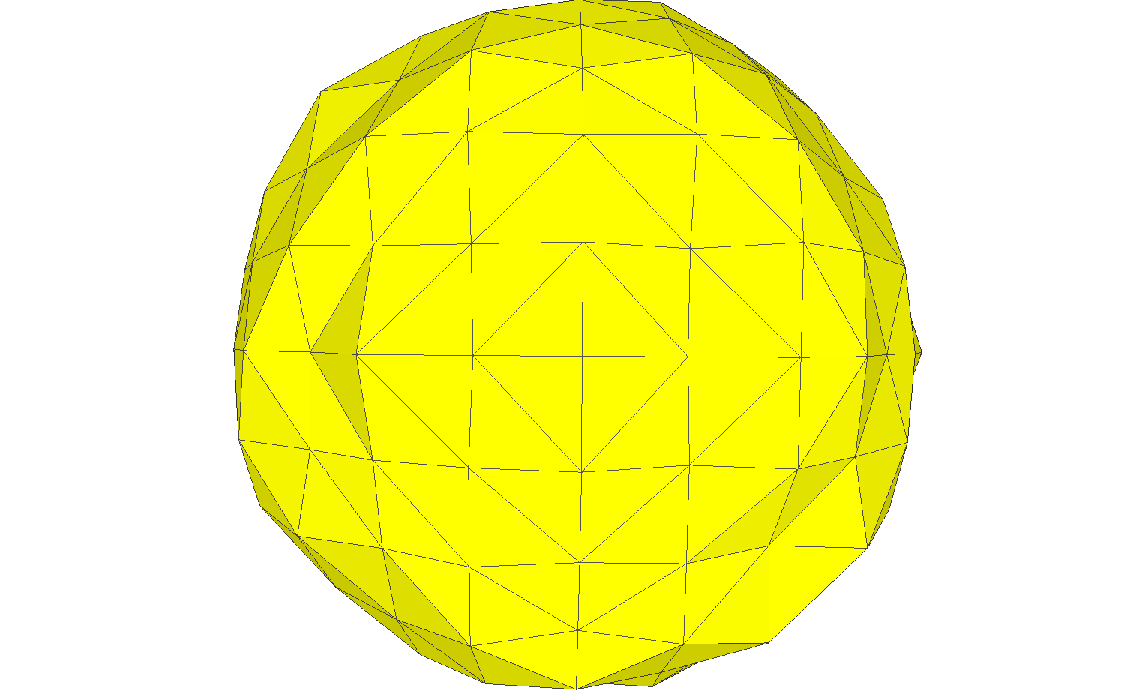}
    \label{fig:subfig1}
}
\subfigure[After 10 steps]{
    \includegraphics[width=0.35\textwidth]{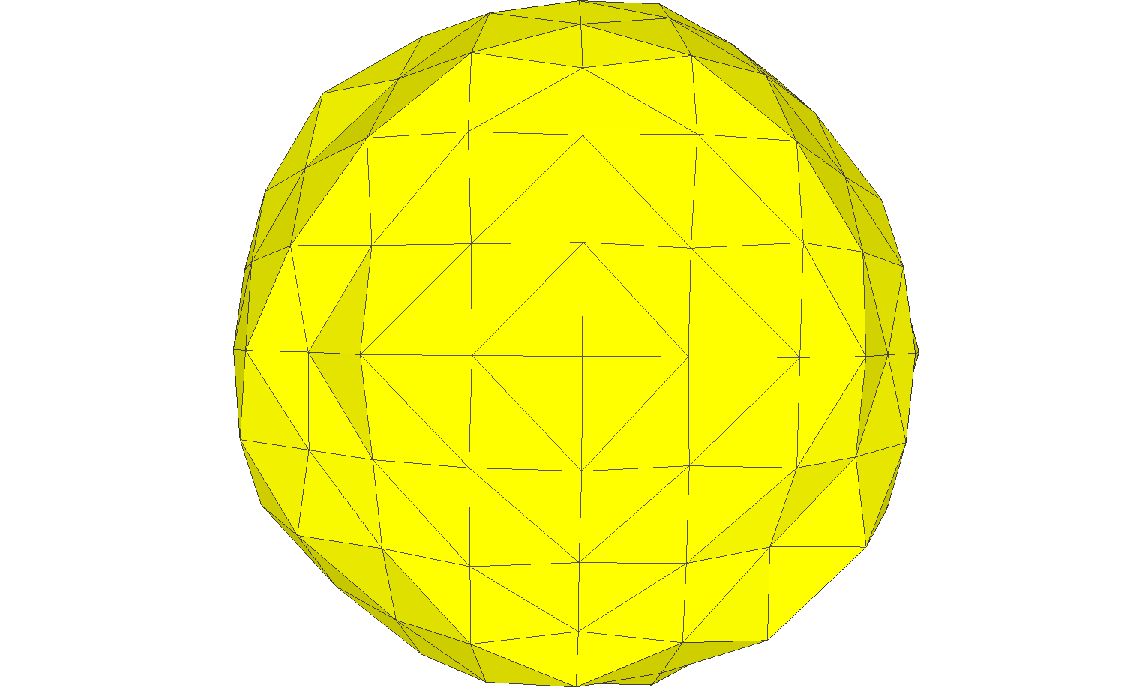}
    \label{fig:subfig2}
}
\subfigure[After 30 steps]{
    \includegraphics[width=0.35\textwidth]{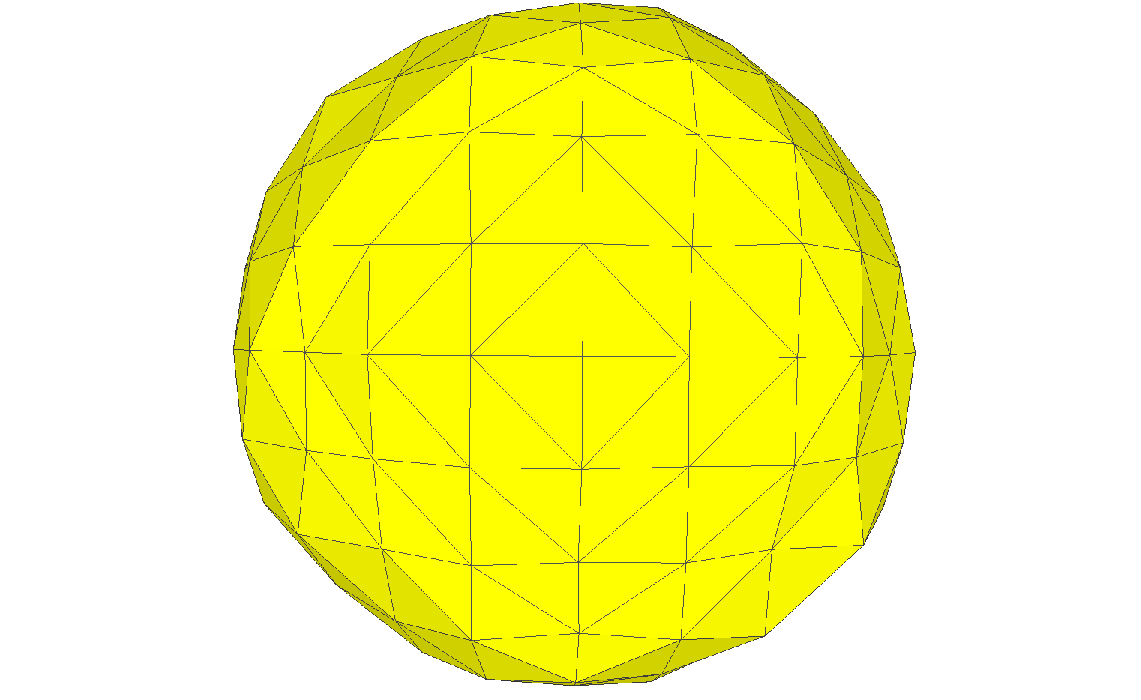}
    \label{fig:subfig3}
}
\subfigure[After 60 steps]{
    \includegraphics[width=0.35\textwidth]{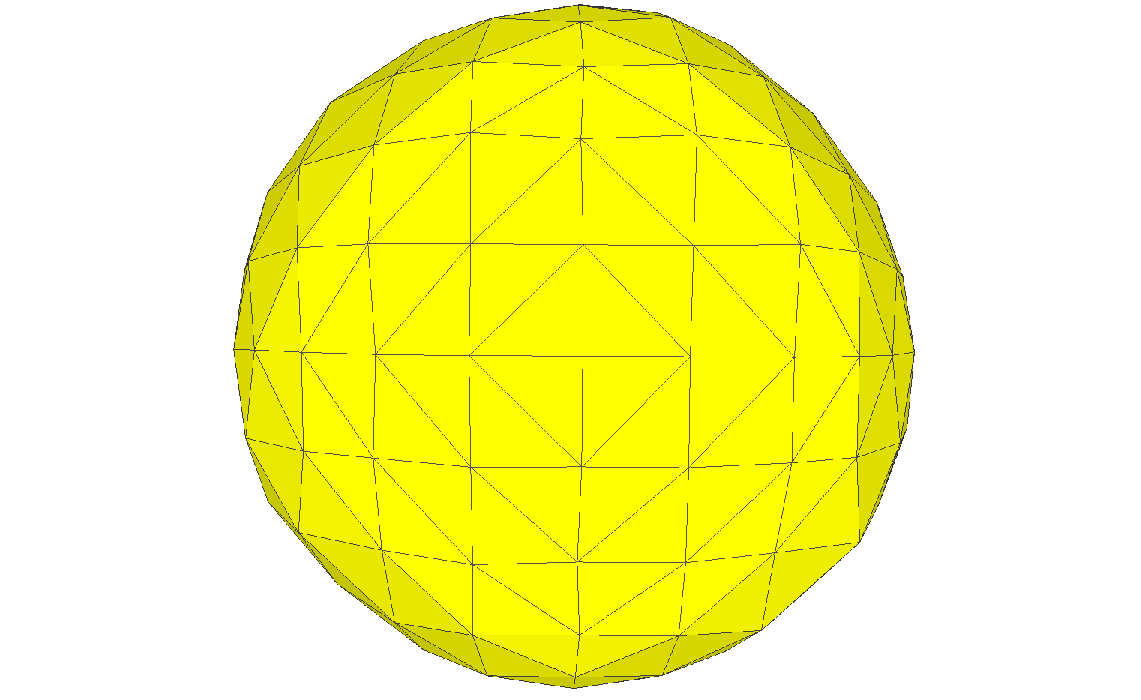}
    \label{fig:subfig4}
}
\caption{Laplacian smoothing on a sphere with cylindrical RBF operator.60 step smoothing was carried out in 0.4 secs including time for computing the laplacian}
\label{laplacesmooth}
\end{figure}

\begin{figure}[ht]
\centering
\subfigure[Histogram of mesh angles]{
    \includegraphics[width=0.45\textwidth]{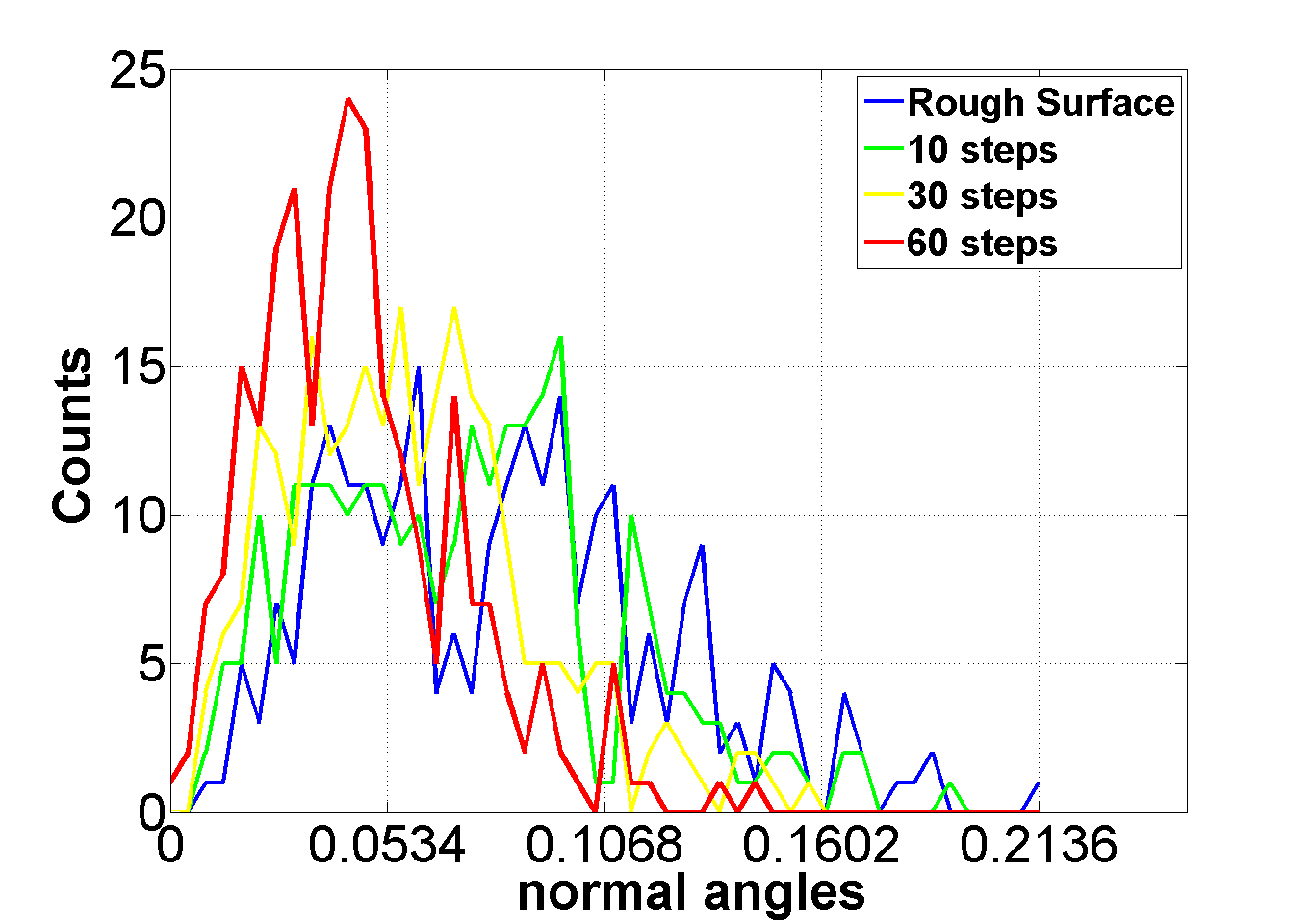}
    \label{histangles}
}
\subfigure[Mean angle trend]{
    \includegraphics[width=0.45\textwidth]{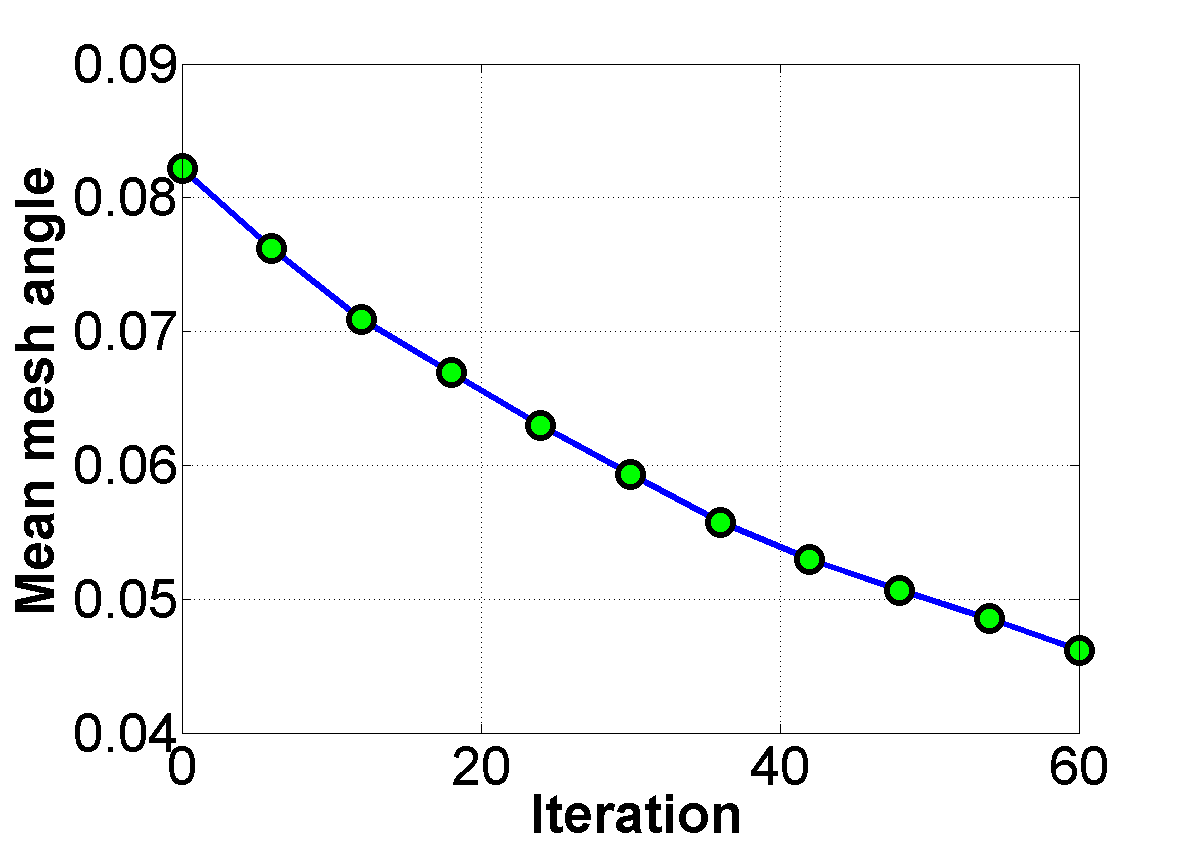}
    \label{meanangles}
}
\subfigure[Mesh smoothing time]{
    \includegraphics[width=0.5\textwidth]{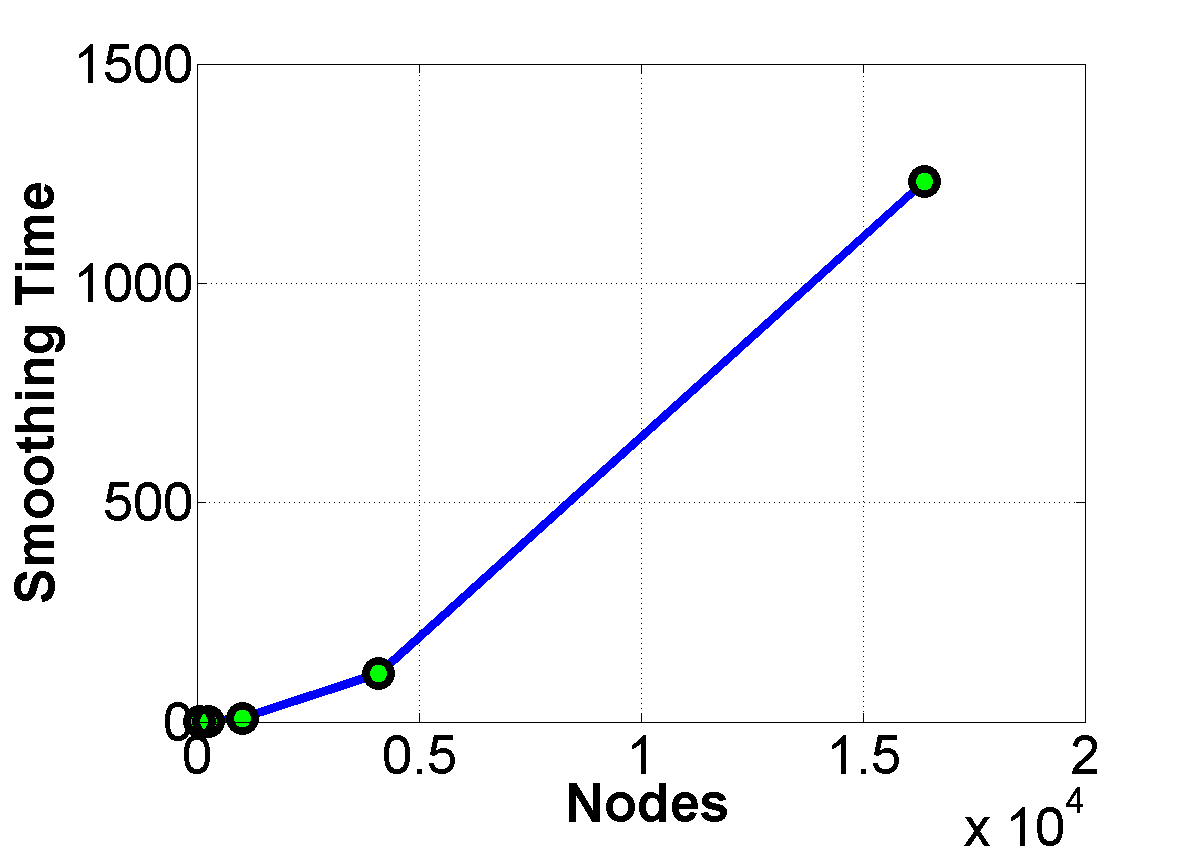}
    \label{smoothtime}
}
\caption{(a) Histogram of angles (in radians) between normals of smoothed sphere (after 10, 30 and 60 steps) and original perfect sphere with 258 points.(b) Trend in mean angles of smoothed sphere. (c) Trend in computational time (in seconds) for 100 steps of CBF smoothing of spheres having nodes ranging between 66-16386.}\label{quantmeshsooth}
\end{figure}

\section{Conclusion}
\label{sect:conclusion}

In this paper we have introduced a new technique for discretizing surface differential operators on closed manifolds. By projecting the radial vector used in standard RBF kernels onto the local tangent plane of the manifold, we produce a modified kernel which is constant in the normal direction to the surface. Our modified kernel, when used to represent functions defined on manifolds, permits the simple replacement of surface differential operators by Cartesian operators. We have demonstrated through numerical experiments the superior performance of CBF in discretizing the Laplace Beltrami operator on the sphere by comparing the first 100 eigen values with the exact values. We also solved the heat equation on the sphere and a molecule and also showed that the discrete operator can be used to effectively smooth noisy surface meshes. Here we worked primarily with the Gaussian kernel. We will consider the effect of other basis functions such as the Matern kernel in future efforts. We will also investigate the approximation of other differential operators such as divergence operator.

\label{sect:bib}
\bibliographystyle{plain}
\bibliography{CRBFbib_easychair}

\end{document}